\documentclass[11pt,a4paper]{article}

 \usepackage[dvipsnames]{color}

     \definecolor{red}{rgb}{0.9,0,0}
     \definecolor{green}{rgb}{0,0.6,0}
     \definecolor{rb}{rgb}{0.6,0,0.2}
     \definecolor{blue}{rgb}{0,0,0.5}
     \definecolor{grey}{rgb}{0.5,0.5,0.5}
\makeatletter
\@addtoreset{equation}{section}
\makeatother

\usepackage{amsmath, amsthm, amssymb}
\usepackage{graphicx}
\newcommand {\eps} {\varepsilon}

\newcommand {\beq} {\begin{equation}}
\newcommand {\eeq} {\end{equation}}
\newcommand {\beqa} {\begin{eqnarray}}
\newcommand {\eeqa} {\end{eqnarray}}
\newcommand {\beqann} {\begin{eqnarray*}}
\newcommand {\eeqann} {\end{eqnarray*}}

\newcommand{\pt}{\partial}

\newtheorem{theorem}{Theorem}[section]

\newtheorem{lemma}[theorem]{Lemma}

\def \Rbb{{\Bbb R}}

\def \i{^{-1}}

\def \half{{1 \over 2}}

\def \rmfor{\ {\rm for} \ }
\def \rmin{\ {\rm in} \ }

\def \rmon{\ {\rm on} \ }

\def \pd{\partial}
\def \<{\langle}
\def \>{\rangle}
\def \nab{\nabla}

\def \LP{\left(}
\def \RP{\right)}

\def \b{\beta}

\def \e{\varepsilon}
\def \d{\delta}
\def \D{\triangle}

\def \g{\gamma}

\def \th{\theta}

\def \x{\xi}

\def \G{\Gamma}

\def \O{\Omega}

\def \y{\eta}

\def \s{\sigma}
\def \t{\tau}

\def \uasS{u_{{\rm as},S}}
\def \bS{\b_S}
\def \uasS{u_{{\rm as},S}}

\def \xb{\bar x}
\def \ev{{\bf e}}
\def \Fc{{\cal F}}
\def \Gc{{\cal G}}
\def \Hc{{\cal H}}

\def \Bt{\tilde B}
\def \vt{\tilde v}

\def \zt{\tilde z}

\def \qt{\tilde q}

\def \vo{\mathring{v}}
\def \vto{\mathring{\vt}}

\title{Some asymptotic expansions for a semilinear reaction-diffusion
problem in a sector%
\thanks{This publication has emanated from research conducted with the
financial support of Science Foundation Ireland under the Basic
Research Grant Programme 2004; Grants 04/BR/M0055 and
04/BR/M0055s1.}}

\author{R. Bruce Kellogg%
\thanks{Department of Mathematics,
University of South Carolina, Columbia, SC 29208, USA;
(rbmjk@windstream.net)}
and
Natalia Kopteva%
\thanks{Department of Mathematics and Statistics,
University of Limerick, Limerick, Ireland (natalia.kopteva@ul.ie)}}

\date{}

\begin{document}

\maketitle


\begin{abstract}
The semilinear reaction-diffusion equation $-\e^2\D u+b(x,u)=0$
with Dirichlet boundary conditions is considered
in a convex unbounded sector.
The diffusion parameter $\eps^2$ is arbitrarily small,
and the ``reduced equation'' $b(x,u_0(x))=0$
may have multiple solutions.
A formal asymptotic expansion for a possible
solution $u$
is constructed that involves boundary and corner layer functions.
For this asymptotic expansion, we establish certain
 inequalities that
are used in \cite{KK0} to construct sharp sub- and
super-solutions and then establish the existence of a solution to
a similar nonlinear elliptic problem in a convex polygon.
\end{abstract}

\section{Introduction}

In this note we consider
the singularly perturbed semilinear reaction-diffusion
boundary-value problem
\begin{subequations}\label{1.1}
\beqa
F u \equiv -\eps^2 \triangle u+ b(x,u) = 0,& \qquad&
x=(x_1,x_2)\in S\subset\mathbb{R}^2,
\label{1.1a}\\
u(x)=g(x),&\qquad&  x\in  \partial S.
\label{1.1b}
\eeqa
\end{subequations}
in a convex sector $S$ with vertex $O$ and sides $\G$ and $\G^-$.
Our purpose in this note
is to establish some asymptotic expansions
and related inequalities for a possible
solution to the problem. These 
are needed in
\cite{KK0} to construct sharp sub- and super-solutions and then
establish the existence of a solution to a similar nonlinear
elliptic problem in a convex polygon.
The proofs involve
lengthy formal calculations, and is the purpose of this paper.

The ``reduced problem'' associated with (\ref{1.1}) is defined by
formally setting $\varepsilon=0$ in (\ref{1.1}a), i.e.
\begin{equation} \label{red}
       b(x,u_0(x))=0 \quad \mbox{for } x\in \bar S.
\end{equation}
It is assumed that \eqref{red} has a smooth solution $u_0$
that is stable in a sense to be described below.
The hypotheses on $b$
are such as to include the possibility of multiple solutions
to (\ref{red}) and therefore to (\ref{1.1}).
Since it may happen that $u_0 \ne g$ on $\pd S $, the solutions
may exhibit boundary layer behavior near $\pd S$.
We shall assume that the function $b$ is smooth
and that $g$ is smooth on each $\G$ and $\G^-$ and continuous at the  vertex $O$.
Furthermore, we assume that $u_0(x)$,  $g(x)$, and, for each fixed $s$,
the function $b(x,s)$,
as well as their derivatives, are bounded as $|x|\rightarrow\infty$.

In addition we make the following assumptions.

\begin{description}
\item[A1]
( {\em stable reduced solution}) There is a number $\g>0$ such that
$$
b_u(x,u_0(x))>\gamma^2>0 \quad \mbox{for all } x\in S.
$$

\item[A2]
({\em boundary condition}) The boundary data $g(x)$ from (\ref{1.1b}) satisfy
$$
\int^{v}_{u_0(x)} \!b(x,s)\,ds>0 \qquad\mbox{for all}\;\;
  v \in \bigl(u_0(x), g(x)\bigr]',\qquad x\in\partial S.
$$
Here the notation $(a,b]'$ is defined to be $(a,b]$ when $a<b$ and
$[b, a)$ when $a>b$, while $(a,b]'=\emptyset$ when $a=b$.

\item[A3]
({\em corner condition})
If $g(O)\neq u_0(O)$, then
$$
\frac{b(O,g(O))}{g(O)-u_0(O)}> 0.
$$

\item[A4]
Only to simplify our presentation, we make a further assumption that
$$
u_0(x) < g(x)\qquad\mbox{for all} \;\;x \in \pd S.
$$
Using A4, we can simplify A3 to
$b( O , g( O ) ) > 0$.

\end{description}

\noindent
Note that if $g(x)\approx u_0(x)$, then A2 follows from A1 combined with (\ref{red}),
while if $g(x)=u_0(x)$ at some point $x\in\pt S$, then A2 does not impose any restriction
on $g$ at this point.
Similarly,  if $g(O)\approx u_0(O)$,
then A3 follows from A1 combined with (\ref{red}),
while if $g(O)=u_0(O)$ at some vertex $O$, then A3 does not impose any restriction
on $g$ at this point.
Assumption A1 is local and permits the construction of multiple solutions
to (\ref{red}) and therefore to (\ref{1.1}).
Assumptions A2 and A3 guarantee existence
of boundary and corner layer ingredients, respectively in an asymptotic expansion
for problem (\ref{1.1}).

The note is organized as follows.
Section~\ref{sec2} defines some boundary layer functions associated
with each side of the sector $S$
and some corner layer functions associated with the vertex of $S$.
The boundary layer
functions are defined as solutions of some ordinary differential equations in a stretched
independent variable. The corner layer functions are solutions of some elliptic partial
differential equations in stretched independent variables.
In Sections~\ref{sec3} and~\ref{sec_beta},
these boundary and corner functions are assembled into
a formal first-order asymptotic expansion
and a perturbed asymptotic expansion, respectively,
and then certain properties of the unperturbed and perturbed asymptotic expansions
are established,
that are used in \cite{KK0}.
The proofs involve much computation, and is the purpose of this note.

\medskip
\noindent
{\it Notation.} Throughout the paper we let $C$, $\bar C$, $c$, $c'$
denote generic positive constants
that may take different values in different formulas, but
are always independent of $\eps$
($\bar C$ is usually used for a sufficiently large constant).
A subscripted $C$ (e.g., $C_1$) denotes a positive constant that is independent of $\eps$
and takes a fixed value. For any two quantities $w_1$ and $w_2$,
the notation $w_1=O(w_2)$ means $|w_1|\le C|w_2|$.



\section{ Boundary and corner layer functions }\label{sec2}
This section defines some boundary layer functions associated with each side of the sector $S$
and some corner layer functions associated with the vertex of $S$.
The boundary layer
functions are defined as solutions of some ordinary differential equations in a stretched
independent variable. The corner layer functions are solutions of some elliptic partial
differential equations in stretched independent variables.
The existence and properties of the
corner layer functions are established in \cite[Section~3]{KK0}.

We use the functions
\beq\label{B_tB}
B(x,t)=b(x,u_0(x)+t),\qquad \Bt(x,t;p)=b(x,u_0(x)+t)-p\,t.
\eeq
The perturbed version $\Bt$ of the function $B$ is used, with $|p|$ sufficiently small,
in the construction of sub-
and super-solutions. In the constructions that follow, a tilde will always denote a perturbed function. The
perturbed functions always depend on the parameter $p$, but we will sometimes not show the explicit
dependence. Thus, we will sometimes write $\Bt(x,t)$ for $\Bt(x,t;p)$.
We need a notation for the derivatives of $\Bt$. For derivatives
with respect to the first argument, we write $\nab_x\Bt$,
$\nab_x^2\Bt$, etc., for the vector, matrix of second derivatives,
etc., with respect to $x$. We write $\Bt_t$, $\Bt_{tt}$, etc., for
derivatives with respect to $t$. Note also that $\Bt(x,0)=0$, so
$\nab_x^k\Bt(x,0)=0$ for $k=1,2,\cdots$, so
\begin{equation}\label{DxB}
|\nab_x^k\Bt(x,t)| \le C|t| \qquad\rmfor k=0,1,2,\cdots.
\end{equation}

We will occasionally use,  for any function $f$, the notations
\begin{equation}\label{fabc}
f\big|^b_a=f(b)-f(a),\qquad f\big|^c_{a;\,b}=f(c)-f(b)-f(a).
\end{equation}
Since $f\big|^{a+b}_{a;\,b}+f(0)=abf''(t)$,
we see that $f(0)=0$ implies
$f\big|^{a+b}_{a;\,b}=O(|ab|)$ and therefore
$f\big|^{a+b+c}_{a;\,b}=O(|c|+|ab|)$.
In view of (\ref{DxB}), we thus have
\beq\label{DxB_abc}
\nab_x^k\Bt(x,\cdot)\Bigr|^{c+a+b}_{a;\,b}=O(|c|+|ab|).
\eeq

We shall now define functions needed to assemble
a first-order
asymptotic expansion and its perturbed version.
The following two subsections deal respectively with a side $\G$ of $S$,
and with the vertex $O$ of $S$.

\subsection{Solution near a side}\label{ssec2_1}

In this subsection we construct boundary layer functions associated with the side $\G$
of $\pd S$. An analogous construction can be made for the side~$\G^-$.
Throughout the subsection,
$\G$ denotes the line that extends the ray~$\G$.
Extend $u_0$ and $b$ to smooth functions, also denoted $u_0$ and $b$,
on $\mathbb{R}^2$ and $\mathbb{R}^2\times \mathbb{R}$, respectively,
so that
(\ref{red}) and A1 hold true for all $x\in\mathbb{R}^2$.
Furthermore,
extend $g$ defined on the ray $\G$ to a smooth function, also denoted
$g$, on the line $\G$,
which satisfies the extended form of A2 and A4 for all  $x \in \G$.

\begin{figure}
\hspace{2.5cm}
\includegraphics{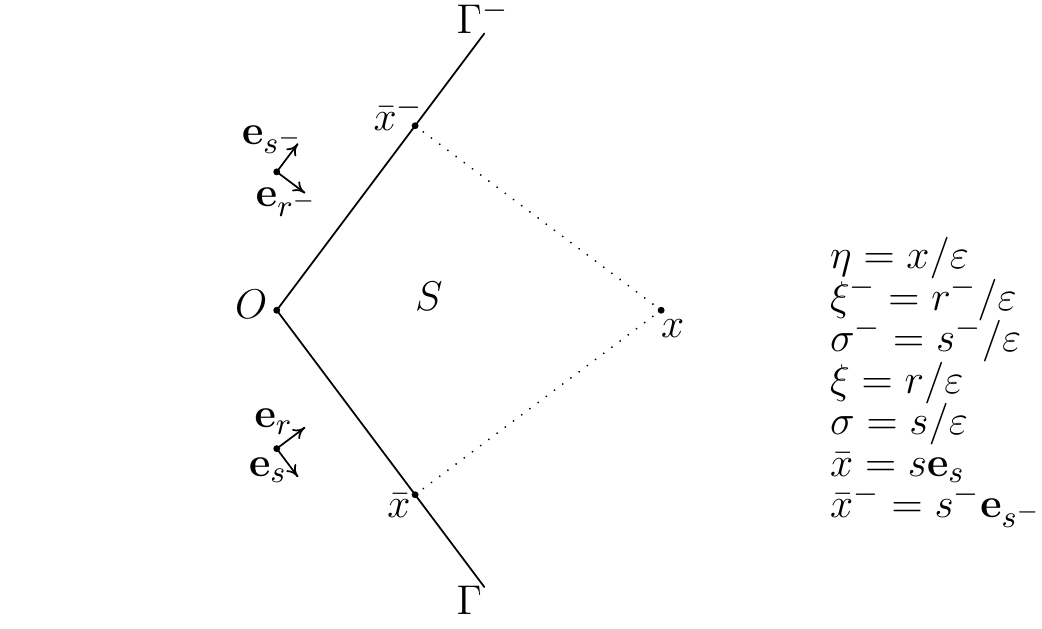}
\caption{Geometry of the sector $S$}\label{fig_fig1}
\end{figure}

Let $\ev_s$ denote the unit vector pointing in the direction of $\G$.
 Let $\ev_r$ be the unit vector perpendicular to $\ev_s$ and oriented
 to point into $S$. Let $s$ denote the signed distance along $\G$ with $s=0$ at $O$
 and $s>0$ on the ray $\G$.
For $x \in \Rbb^2$ write $x=O+s\ev_s+r\ev_r$.
Then $\xb=O+s\ev_s$ is the point on $\G$ which is closest to $x$ and $r$ is
the signed distance from $\xb$ to $x$, with $r>0$ if $x \in \O$.
($\ev_s$, $\ev_r$, $x$ and $\bar x$ are shown
in Figure~\ref{fig_fig1}).

Let $\vt_0(\xi,s;p)$ be the solution
to the nonlinear autonomous two point boundary value problem
\begin{subequations}\label{v0}
\beqa
\label{v0a}
&\displaystyle-\frac{\pt^2 \vt_0}{\pt \xi^2}+\Bt(\bar x,\vt_0;p)=0,\\
\displaystyle\label{v0b}
&\vt_0(0,s;p)=g(\bar x)-u_0(\bar x), \qquad \vt_0(\infty,s;p)=0.
\eeqa
\end{subequations}
The geometric meaning of the variable $\x$ is given by the formula
$\x=r/\e$. The variables $p$ and $s$ appear as parameters in the
problem \eqref{v0}.
The parameter $p$ satisfies $|p|<\g^2$ and in
general will be close to zero. We sometimes omit the explicit
dependence of $\vt_0$ on $p$ and write $\vt_0(\x,s)=\vt_0(\x,s;p)$.

We set $v_0(\x,s)=\vt_0(\x,s;0)$. The function  $v_0$ appears in the asymptotic expansion
of the solution near the side $\G$.
With $v_0$ defined, we define a function $v_1(\x,s)$ to be the solution to the linear two point
boundary value problem
\begin{subequations}\label{v1}
\beqa
&\displaystyle
-\frac{\pt^2v_1}{\pt \xi^2}+v_1 B_t(\bar x,v_0) =-\x\,\ev_r\cdot\nab_xB(\xb,v_0), \label{v1a}\\
\displaystyle
& v_1(0,s)=v_1(\infty,s)=0.\label{v1b}
\eeqa
\end{subequations}
Note that $v_1$ is not a perturbed function as it does not depend on $p$.
We also define
\begin{equation}\label{v}
\begin{array}{c}
\vto_0(\x;p)=\vt_0(\x,0;p),\qquad \vo_0(\x) = v_0(\x,0),\qquad \vo_1(\x)=v_1(\x,0), \\
\vt = \vt_0+\e v_1,\quad\; v=v_0+\e v_1, \quad\; \vto=\vto_0+\e \vo_1,
\quad\; \vo=\vo_0+\e \vo_1. \\
\end{array}
\end{equation}
In our notation, a small circle above a function name indicates that in the argument of the function
we have set $s=0$.

For the solvability and properties of problems \eqref{v0} and \eqref{v1} we
cite a result from \cite[Lemma~2.1]{KK0}.

\begin{lemma}\label{l-v}
There is $p_0\in(0,\gamma^2)$
such that for all $|p|\le p_0$
%
there exist functions $\vt_0$ and $v_1$
that satisfy \eqref{v0},\,\eqref{v1}.
For the function $\vt_0=\vt_0(\xi,s;p)$ we have
\beq\label{vt_monotone}
\vt_0\ge 0, \qquad\quad \frac{\partial \vt_0}{\partial p}\ge 0.
\eeq
Furthermore, for any $k \ge 0$ and arbitrarily small but fixed $\d$,
there is a $C>0$ such that for
$0\le \xi<\infty $,
$s \in \Rbb$ and $k=0,1,\cdots$,
\begin{equation*}
\Bigl|\frac{\partial^k  \vt_0}{\partial \xi^k} \Bigr|
+\Bigl|\frac{\partial^k \vt_0}{\partial s^k}\Bigr|
+\Bigl|\frac{\partial^k  v_1}{\partial \xi^k} \Bigr|
+\Bigl|\frac{\partial^k v_1}{\partial s^k}\Bigr|
+\Bigl|\frac{\partial \vt_0}{\partial p}\Bigr|
+\Bigl|\frac{\partial^2 \vt_0}{\partial p\,\partial s}\Bigr|
 \le C  e^{-(\gamma-\sqrt{|p|}-\delta)\xi}.
\end{equation*}
\end{lemma}

For later purposes we shall now obtain an estimate for $\vt_0-v_0$.

\begin{lemma}\label{vt0-v0}
We have, for $|p|$ sufficiently small,
\begin{equation}\label{Dw0ineq}
-\e^2\D (\vt_0-v_0) =-B(x,\cdot)\Bigl|^{\vt}_{v}
+p v_0+O(\e^2+p^2).
\end{equation}
\end{lemma}

\begin{proof}
It follows from Lemma~\ref{l-v}
that
$|\vt_0-v_0| \le Cpe^{-c\x}$ and $|\frac{\pt^2}{\pt s^2}\vt_0|\le C$.
The latter estimate yields
$\e^2\D (\vt_0-v_0)={\textstyle\frac{\pt^2}{\pt\xi^2}}(\vt_0-v_0)+O(\e^2)$.
Furthermore, invoking (\ref{v0a}),\,(\ref{B_tB}) and the estimate for $\vt_0-v_0$, we get
$$
-{\textstyle\frac{\pt^2}{\pt\xi^2}}(\vt_0-v_0)
=-B(\xb,\cdot)\Bigl|^{\vt_0}_{v_0}+p\vt_0
=-\Gc(\xb,0)+pv_0+O(p^2),
$$
where we use the auxiliary function $\Gc(x,t)=B(x,\cdot)\Bigl|^{\vt_0+t}_{v_0+t}$.
Noting that
$B(x,\cdot)\Bigl|^{\vt}_{v}
=B(x,\cdot)\Bigl|^{\vt_0+\eps v_1}_{v_0+\eps v_1}=\Gc(x,\eps v_1)$,
it remains to establish the estimate $\Gc(x,\eps v_1)-\Gc(\xb,0)=O(\e^2+p^2)$.
Indeed, we have
\beq\label{last}
\Gc(x,\eps v_1)-\Gc(\xb,0)
=(x-\xb)\cdot\nab_x\Gc^*+\e v_1\Gc^*_t
=O(\e^2+p^2),
\eeq
where the asterisks indicate that the derivatives are evaluated at an
intermediate point. In the last step in (\ref{last}) we combined
 $|\Gc_t^*|=|\Gc_t(x^*,t^*)|=|B_t(x^*,\vt_0+t^*)-B_t(x^*,v_0+t^*)|
\le C|\vt_0-v_0|$,
and a similar estimate
$|\nab_x\Gc^*| \le C|\vt_0-v_0|$
with $|x-\xb|=\e\x$, $|v_1| \le C$ and
$|\vt_0-v_0| \le Cpe^{-c\x}$.
\end{proof}

\subsection{Solution near a vertex}\label{ssec2_2}

In this subsection we construct corner layer functions
associated with the vertex $O$.
Some notation is required for the constructions.
Let $s$ denote the distance along $\G$,
measured from $O$, and let $r$ denote the perpendicular distance to a point $x \in S$.
Thus, $x \to (s,r)$ is a linear orthogonal map.
We also let $\ev_s$ and $\ev_r$ denote the unit vectors along $\G$ and orthogonal to $\G$
respectively, so $x=r\ev_r+s\ev_s$.
We denote by $\xb=s\ev_s$ the point of $\G$ that is
closest to $x$.
In a similar manner, we define variables
$(s^-,r^-)$,
so $x=r^-\ev_{r^-}+s^-\ev_{s^-}$,
and $\xb^-=s^-\ev_{s^-}$ associated with the side $\G^-$.
The variable $s^-$ denotes the distance along $\G^-$, measured from $O$.
We will also need stretched variables.
We set $\y=x/\e$, $\x=r/\e$, $\s=s/\e$, $\x^-=r^-/\e$, $\s^-=s^-/\e$. These variables
are shown in Figure~1.

Using these notations, Section~\ref{ssec2_1}
gives functions $\vt_0(\x,s;p)$ and $v_1(\x,s)$ associated
with the side $\G$ and functions $\vt_0^-(\x^-\!\!,s^-\!;p)$ and $v_1^-(\x^-\!\!,s^-)$
associated
with the side $\G^-$.
We also recall the notations in \eqref{v} and use corresponding notations for the
side $\G^-$. The function $\vt$ matches the
disparity between the boundary conditions of (\ref{1.1}b) and the
value of $u_0$ on $\G$, but leaves a rapidly decaying boundary
value on $\G^-$. The function $\vt^-$ has a similar behavior, with a rapidly
decaying boundary value on $\G$. To
deal with these rapidly decaying boundary values we construct
functions $\zt_{0}(\y;p)$ and $z_{1}(\y)$, defined in terms of the
stretched variable $\y$.

The function $\zt_{0}$ is defined to be
a bounded solution of the
autonomous nonlinear elliptic boundary value problem
%
\begin{equation}\begin{array}{rcl}\label{z0}
-\D_\y \zt_{0}+\Bt(O,\zt_{0};p)=0&\quad& \rmin S, \\
\zt_{0}=A:=g(O)-u_0(O) &&\rmon \pd S.
\end{array}\end{equation}
%
Here we have $A>0$, by our
assumption A4 at the point $O$.
We also set $z_{0}(\y)=\zt_{0}(\y;0)$.
The existence and properties of $z_{0}$ are given in the following theorem;
see \cite[Theorem~2.2]{KK0}.

\def \dist{{\rm dist}}
\begin{theorem}\label{T-z0}
There is a positive constant $p^*$ such that if $|p| \le p^*$,
the problem \eqref{z0} has, for each $p$, a solution $\zt_{0}$ which satisfies
$\zt_{0}\le A$ and
\begin{equation}\label{T-z0-a}
0<\max\{\vto_0,\vto_0^-\} \le \zt_{0}(\y;p) \le
\max\{\vto_0,\vto_0^-\}+C|\eta|^{-1},
\end{equation}
and which is an increasing function of $p$.
Also, $|\nab\zt_0| $ is bounded in $S$.
Finally there is a constant $C>0$ such that
\begin{equation}\label{T-z0-b}
\zt_{0}(\y) \le C\LP e^{-\g\x}+e^{-\g\x^-} \RP .
\end{equation}
\end{theorem}

We also consider a function $z_1(\y)$ which satisfies the linear
elliptic boundary value problem
\beq\label{z1}
\begin{array}{c}
-\D_\y z_{1}+z_{1}B_t(O,z_{0}) =-\eta\cdot\nabla_x B(O,z_{0}) \;\;\rmin S, 
\\
z_{1}=\s\,
{\textstyle\frac{\pd}{\pd s}}(g-u_0)\bigr|_{x=O}
\;\rmon \G,\quad\quad
z_{1}=\s^-
{\textstyle\frac{\pd}{\,\pd s^-\!}}(g^--u_0)\bigr|_{x=O}
\;\rmon \G^-, 
\end{array}
\eeq

The functions $\zt_0$ and $z_1$
form a correction $\zt_0+\eps z_1$ to the reduced solution $u_0$
in close proximity of the vertex $O$.
To extend it further away from $O$, the corrections
$\vt_0+\e v_1$ and $\vt^-_0+\e v_1^-$ to $u_0$ near the sides $\G$ and $\G^-$ are to be invoked
as follows.
We use the corner functions $\zt_0$ and $z_1$
together with the boundary functions $\vt_0$, $v_1$, $\vt^-_0$, $v_1^-$
to define a related pair of corner functions $\tilde q_0$ and $q_1$,
which, rather than $\zt_0$ and $z_1$,
will appear in a formal asymptotic expansion of the
solution of \eqref{1.1} in the entire $S$;
see Sections~\ref{sec3},\,\ref{sec_beta} below.

We shall use the following notation.
Pick a point $\y \in S$. Having chosen $\y$, the
formulas
\begin{equation}\label{yrs}
\y=\x\ev_r+\s\ev_s = \x^-\ev_{r^-}+\s^-\ev_{s^-}
\end{equation}
determine numbers $\x,\s,\x^-,\s^-$; see Figure~1.
With this notation, and using
the functions  $\zt_{0}$, $z_{1}$
and $\vto_0$, $\vto_0^-$, $\vo_1$, $\vo^-_1$ of \eqref{v0},\,\eqref{v1},\eqref{v},
we  define
\begin{subequations}\label{qdef}
\begin{align}
\qt_{0}(\y;p) &= \zt_{0}(\y;p)-\vto_0(\x;p)-\vto_0^-(\x^-;p),
\label{qdefa} 
\\
q_{1}(\y) &= z_{1}(\y)
-[\vo_1(\x)+ \s\vo_{0,s}(\x)]
-[\vo_1^-(\x^-) + \s^- \vo^-_{0,s^-}(\x^-)],
\label{qdefb}
\end{align}
and furthermore,
\beq\label{qdefc}
\!
\qt(\y;p) = \qt_{0}(\y;p)+\e q_{1}(\y), 
\;\;\;
q_{0}(\y)= \qt_{0}(\y;0),\;\;\;  q(\y)=q_{0}(\y)+\e q_{1}(\y).%
\eeq
In these formulas, following the notational conventions of \eqref{v},
we mean
\beq\label{qdefd}
\vo_{0,s}(\x)={\textstyle\frac{\pd}{\pd s}}v_0(\x,s)\bigr|_{s=0},
\qquad
\vo^-_{0,s^-}={\textstyle\frac{\pd}{\,\pd s^-\!}}v^-_0(\x^-\!\!,s^-)\bigr|_{s^-=0}.
\eeq
\end{subequations}
Under this notation, the boundary conditions in (\ref{z1}) become
\beq\label{z1_bc}
z_{1}=\s\,\vo_{0,s}\quad  \rmon \G,\qquad\quad
z_{1}=\s^- \vo^-_{0,s^-} \quad \rmon \G^-.
\eeq

From the above formulas,
noting that
$\D_\y\qt_0=\D_\y\zt_0-\frac{\pt^2}{\pt\xi^2}\vto_{0}-\frac{\pt^2}{(\pt\xi^-)^2}\vto^-_{0}$,
and using \eqref{v0},\,\eqref{z0},
we derive a nonlinear boundary value problem satisfied by $\qt_0$:
\begin{subequations}\label{qt0}
\begin{align}
\D_\y\qt_0 &=  \Bt(O,\qt_0+\vto_0+\vto^-_0)-\Bt(O,\vto_0)- \Bt(O,\vto_0^-) , \label{qt0a} \\
\qt_0 &= -\vto_0^- \rmon \G,\qquad\quad \qt_0 = -\vto_0 \rmon \G^-. \label{qt0b}
\end{align}
\end{subequations}
Similarly (see Lemma~\ref{bvps-q} below for details),
using \eqref{v0}, \eqref{v1} and \eqref{z1},
we formally derive a linear boundary value problem satisfied by $q_1$:
\begin{subequations}\label{q1bvp}
\begin{align}
-\D_\y q_1&+q_{1}B_t(O,z_{0}) =-\y\cdot\nab_x B(O,\cdot)\Bigr|_{\vo_0;\ \vo_0^-}^{z_{0}} \notag\\
   &-(\vo_1+\s\vo_{0,s}) \  B_t(O,\cdot)\Bigr|_{\vo_0}^{z_{0}}
 -(\vo_1^- +\s^-\vo^-_{0,s^-})\  B_t(O,\cdot)\Bigr|_{\vo_0^-}^{z_{0}},  \label{q1bvpa}\\
 q_1 = &-(\vo_1^-+\s^-\vo_{0,s^-}^-) \rmon \G,
 \qquad
q_1 = -(\vo_1+\s\vo_{0,s}) \rmon \G^-. \label{q1bvpb}
\end{align}
\end{subequations}
where we used the notation~\eqref{fabc}.
%
%
%
%
Finally, by formally differentiating relation \eqref{qdefa} and problem \eqref{z0}
(or the equivalent problem \eqref{qt0}) with respect to $p$
and invoking (\ref{B_tB}),
we formally derive a boundary value problem that is satisfied by $\qt_{0,p}$:
\beq\label{qt0pbvp}
\begin{array}{c}
-\D_\y\qt_{0,p}+\qt_{0,p}\Bt_t(O,\zt_0)
  = \qt_0
  -\vto_{0,p}\ \Bt_t(O,\cdot)\Bigr|^{\zt_0}_{\vto_0}
  -\vto^-_{0,p}\ \Bt_t(O,\cdot)\Bigr|^{\zt_0}_{\vto^-_0}\ ,
   \\
\qt_{0,p} = -\vto^-_{0,p} \rmon \G,\qquad\quad \qt_{0,p} = -\vto_{0,p} \rmon \G^-.
\end{array}
\eeq

It is shown in \cite[Lemmas~3.6,\ 3.16,\ 3.17]{KK0} that the functions
$\qt_0$, $q_1$ and $\qt_{0,p}$ exist and
are exponentially decaying in $S$, i.e.
there are constants $C_1$ and $c_1$ such that
\begin{equation}\label{q_decay}
|\qt_0|+|q_1|+|\qt_{0,p}|
\le C_1e^{-c_1|\y|} \qquad\rmin S.
\end{equation}
In view of \eqref{qdefb}, the existence of $q_1$ immediately implies existence of $z_1$.
Similarly, having proved the existence
of the solution to \eqref{qt0pbvp}, an integration is used to show that this solution is in fact
the derivative of $\qt_0$ with respect to $p$.

We now derive the boundary value problem satisfied by the function $q_1$.

\begin{lemma}\label{bvps-q}
The function $q_1$ defined by (\ref{qdefb})
satisfies problem (\ref{q1bvp}).
\end{lemma}

\begin{proof}
To prove \eqref{q1bvpa}, note that
\beq\label{r1}
\D_\y q_1 =  \D_\y z_1
-(\D_\y \vo_1+ \s\,\D_\y\vo_{0,s})
-(\D_\y \vo^-_1 + \s^-\D_\y \vo^-_{0,s^-}).
\eeq
Next, using (\ref{qdefd}) and then (\ref{v0a}), we calculate
$$
\D_\y\vo_{0,s}={\textstyle\frac{\partial}{\partial s}} v_{0,\xi\xi}\bigr|_{s=0}
={\textstyle\frac{\partial}{\partial s}} B(s\ev_s,v_0)\bigr|_{s=0}
=\ev_s\cdot\nab_xB(O,\vo_0)+ \vo_{0,s}B_t(O,\vo_0).
$$
Combining this with \eqref{v1a}, yields
\beqa
\D_\y \vo_1+ \s\,\D_\y\vo_{0,s}&=&
[\vo_1B_t(O,\vo_0)+\x\ev_r\cdot \nab_xB(O,\vo_0)]\nonumber\\\nonumber
&&\;\;{}+\s[\ev_s\cdot\nab_xB(O,\vo_0)+ \vo_{0,s}B_t(O,\vo_0)]\\
&=&(\vo_1+ \s\vo_{0,s})B_t(O,\vo_0)
+\eta\cdot \nab_xB(O,\vo_0),
\label{r2}
\eeqa
where we used $\x\ev_r+\s\ev_s=\y$ from (\ref{yrs}).
Similarly, one gets
\beq\label{r3}
\D_\y \vo^-_1+ \s^-\D_\y\vo^-_{0,s^-}=(\vo^-_1+ \s^-\vo^-_{0,s^-})B_t(O,\vo^-_0)
+\eta\cdot \nab_xB(O,\vo^-_0).
\eeq
Recalling that, by (\ref{z1}), we have
$\D_\y z_1=z_1  B_t(O,z_0 ) + \y \cdot \nab_x B(O,z_0)$, where
in the right-hand side $z_1$ is replaced by
$q_1+(\vo_1+  \s\vo_{0,s})+(\vo^-_1  + \s^- \vo^-_{0,s^-})$,
and combining this with (\ref{r1}),\,(\ref{r2}),\,(\ref{r3}), yields \eqref{q1bvpa}.
Finally, noting that $\vo_1=0$ on $\G$ and $\vo_1^-=0$ on $\G^-$,
and then comparing (\ref{qdefb}) and (\ref{z1_bc}), we immediately get \eqref{q1bvpb}.
\end{proof}

\section{Asymptotic expansion}\label{sec3}

In  Section~\ref{ssec2_1} we have defined boundary layer functions
$\vt=\vt_{0}+\eps v_{1}$  and $\vt^-=\vt^-_{0}+\eps v_{1}^-$ associated, respectively,
with
the sides $\G$ and $\G^-$ of $S$,
and in Section~\ref{ssec2_2} we have defined corner layer functions
$\qt=\qt_{0}+\e q_{1}$ associated with the vertex $O$ on $S$.
In the present and next sections,
these functions are used to assemble a formal first-order asymptotic expansion
and then a perturbed asymptotic expansion
for the problem \eqref{1.1}.
We establish certain properties of
the unperturbed and perturbed asymptotic expansions
that are used in \cite{KK0}.
The proofs involve
lengthy formal calculations, and is the purpose of this paper.

The  asymptotic expansion $\uasS$ is defined as follows:
\beqa
\uasS(x) & =&u_0(x)+v(\x,s)+v^-(\x^-,s^-)+q(\y). \label{uas}
\eeqa
The next lemma shows that
the differential equation applied to this
asymptotic expansion is $O(\e^2)$.

\begin{lemma}\label{Fuas}
For the asymptotic expansion $\uasS(x)$ of (\ref{uas}) one has
\begin{subequations}\label{resid}
\begin{align}
F\uasS&=O(\e^2),  \label{resida}\\
\uasS(x)&=g(x)+O(\eps^2) \rmfor x \in \pd S. \label{residb}
\end{align}
\end{subequations}
\end{lemma}

\begin{proof}
(i) We start by establishing
\beq\label{B(x,v)}
-\eps^2\D v+B(x,v)=O(\eps^2),
\qquad
-\eps^2\D v^-+B(x,v^-)=O(\eps^2).
\eeq
The first bound here is obtained  estimating $B(x,v)$ as follows.
Fix $\xi$ and $s$; then $B(x,v)$
is a function of $\eps$, i.e.
$B(x,v)=\Gc(\e)$, where
\begin{equation*}
\begin{aligned}
\Gc(\e)=B\big(\xb+\e\x\ev_r,v_0 +\e v_1\big)
\quad\mbox{with}
\; \xb=s\ev_s,\;v_0=v_0(\x,s),\;v_1=v_1(\x,s).
\end{aligned}
\end{equation*}
Expand $\Gc$ in a Taylor series around $\e=0$ to obtain
\beq\label{G_eps}
\Gc(\e)=\Gc(0)+\e \Gc'(0)+{\textstyle\half}\e^2\Gc''(\e^*)
\eeq
with $0<\e^*<\e$. A calculation shows that
\begin{equation*}
\begin{aligned}
\Gc'(\e)&=\x \ev_r\cdot\nab_x
   B(\xb+\e\x\ev_r,v_0 +\e v_1) + v_1B_t(\xb+\e\x\ev_r,v_0 +\e v_1), \\
\Gc''(\e)&=\x^2 \ev_r^T[\nab_x^2B(\xb+\e\x\ev_r,v_0+\e v_1)]\ev_r \\
 &\quad + 2v_1\x \ev_r\cdot\nab_x B_t(\xb+\e\x\ev_r,v_0+\e v_1)
  + v_1^2B_{tt}(\xb+\e\x\ev_r,v_0 +\e v_1). \\
\end{aligned}
\end{equation*}
Hence
\begin{equation}\label{G_eps_0}
\begin{aligned}
\Gc(0) &= B(\xb,v_0)={\textstyle\frac{\pt^2}{\pt\x^2}}v_{0}, \\
\Gc'(0)&=\x \ev_r\cdot\nab_x B(\xb,v_0 )  + v_1B_t(\xb,v_0 )
={\textstyle\frac{\pt^2}{\pt\x^2}}v_{1},
\end{aligned}
\end{equation}
where we also used \eqref{v0a} and \eqref{v1a}.
Applying \eqref{DxB}, we get
\begin{equation*}
|\nab_x^2B(\xb+\e\x\ev_r,v_0+\e v_1)| \le C|v_0+\e v_1|.
\end{equation*}
By Lemma~\ref{l-v}, $v_0$ and $v_1$ are exponentially decaying in $\x$,
which yields
$\x^2|v_0+\e v_1| \le C$. Hence the first term in the formula for
$\Gc''(\e)$ is bounded. The other 2 terms are bounded for a
similar reason, so $|\Gc''(\e)| \le C$,
where $C$ is independent of $\xi$ and $s$.
Combining this with (\ref{G_eps}) and (\ref{G_eps_0}), yields
\begin{equation*}
B(x,v)=\Gc(\e)
= {\textstyle\frac{\pt^2}{\pt\x^2}}(v_{0} + \e v_1)  +O(\e^2)
={\textstyle\frac{\pt^2}{\pt\x^2}}v  +O(\e^2).
\end{equation*}
As
$\eps^2\D v={\textstyle\frac{\pt^2}{\pt\x^2}}v+O(\eps^2)$,
the first bound in (\ref{B(x,v)}) is established.
The second bound is obtained similarly.


(ii) To show \eqref{resida}, we calculate
\beqa
F\uasS &=& -\e^2\D[u_0+v+v^-+q]+b(x,u_0+v+v^-+q) \nonumber\\\label{ineqsb1}
  &=& B(x,\cdot)\Big|^{v+v^-+q}_{v;\ v^-}
  -\D_\y q  + O(\e^2),
\eeqa
where we used (\ref{B_tB}) and \eqref{B(x,v)}, and also
the notation (\ref{fabc}).

Fix a point $\y \in S$;
then the first term in \eqref{ineqsb1} is a function of $\eps$,
which we denote $\Fc(\e)$. To be more precise, having chosen $\y$, the
formulas \eqref{yrs}
determine fixed numbers $\x,\s,\x^-,\s^-$.
With the understanding
that
\begin{subequations}\label{underst}
\begin{equation}\label{underst_a}
\begin{aligned}
&v =v(\e)= v_0(\x,\e\s)+\e v_1(\x,\e\s), \\
&v^- =v^-(\e)= v^-_0(\x^-,\e\s^-)+\e v^-_1(\x^-,\e\s^-), \\
&q =q(\e)= q_0(\y)+\e q_1(\y),
\end{aligned}
\end{equation}
define a function $\Fc(\e)$ by
\begin{equation}\label{underst_b}
\Fc(\e)=B(\e\y,\cdot)\Bigr|_{v;\ v^-}^{v+v^-+q}.
\end{equation}
\end{subequations}
In view of \eqref{ineqsb1}, to prove \eqref{resida}, we need to show that
\begin{equation*}
\begin{aligned}
\Fc(\e)-\D_\y [q_0+\e q_1] = O(\e^2) \qquad\rmin S.
\end{aligned}
\end{equation*}
Thus  we must show that there is a number $C$, independent
of $\y$, such that
\begin{subequations}\label{ineqsb2}
\begin{align}
\Fc(0) &= \D_\y q_0,\label{ineqsb2a} \\
\Fc'(0) &= \D_\y q_1,\label{ineqsb2b} \\
|\Fc''(\e)| &\le C. \label{ineqsb2c}
\end{align}
\end{subequations}
From the definition (\ref{underst}) of $\Fc$, we have
$v\bigr|_{\eps=0}=\vo_0$, $v^-\bigr|_{\eps=0}=\vo^-_0$,
and hence
\begin{equation*}
\begin{aligned}
\Fc(0)=B(O,\vo_0+\vo_0^-+q_0)-B(O,\vo_0)-B(O,\vo^-_0),
\end{aligned}
\end{equation*}
so \eqref{qt0a} gives \eqref{ineqsb2a}.


A calculation using (\ref{underst_a}) yields
\begin{equation}\label{dv_de}
\frac{dv}{d\e}=v_1 +\s v_{0,s}+\e\s v_{1,s}\,,
\qquad
\frac{d^2v}{d\e^2}=2\s v_{1,s} +\s^2 v_{0,ss}+\e\s^2 v_{1,ss}\,,
\end{equation}
similar relations for $v^-$,
and also $\frac{dq}{d\e}=q_1$, $\frac{d^2q}{d\e^2}=0$.
As $|\sigma|\le |\eta|$ and $|\sigma^-|\le|\eta|$, invoking
Lemma~\ref{l-v} and (\ref{q_decay}),
for $k=0,1,2$ we get
\begin{equation}\label{dv_de_est}
\bigl|\frac{d^kv}{d\e^k}\bigr|\le C(1+|\eta|^k)e^{-c\xi},\;
\bigl|\frac{d^kv^-}{d\e^k}\bigr|\le C(1+|\eta|^k)e^{-c\xi^-},\;
\bigl|\frac{d^kq}{d\e^k}\bigr|\le Ce^{-c|\eta|}.
\end{equation}


We now calculate
\begin{equation*}
\begin{aligned}
\Fc'(\e) &=\y\cdot\nab_x B(\e\y,\cdot)\Bigr|_{v;\ v^-}^{v+v^-+q}\\
&+\frac{dq}{d\e}B_t(\e\y,\cdot)\Bigr|_{v+v^-+q}
+\frac{dv}{d\e}B_t(\e\y,\cdot)\Bigr|_{v}^{v+v^-+q}
+\frac{dv^-\!\!\!}{d\e}B_t(\e\y,\cdot)\Bigr|_{v^-}^{v+v^-+q}.
\end{aligned}
\end{equation*}
Hence, using the first relation in (\ref{dv_de}) and its analogue for $v^-$, we get
\begin{equation*}
\begin{aligned}
\Fc'(0)&=\y\cdot\nab_x B(O,\cdot)\Bigr|_{\vo_0;\ \vo_0^-}^{\vo_0+\vo_0^-+q_0}\\
&+q_1B_t(O,\vo_0+\vo_0^-+q_0)\\
&+(\vo_1+\s\vo_{0,s})B_t(O,\cdot)\Bigr|_{\vo_0}^{\vo_0+\vo_0^-+q_0}
+(\vo^-_1+\s^-\vo^-_{0,s})B_t(O,\cdot)\Bigr|_{\vo_0^-}^{\vo_0+\vo_0^-+q_0}.
\end{aligned}
\end{equation*}
Recalling that $\vo_0+\vo_0^-+q_0=z_0$ and
inspecting \eqref{q1bvpa} we see that (\ref{ineqsb2}b) holds.


A formula for the quantity $\Fc''(\e)$ is
obtained by a lengthy but straightforward computation, which gives
\begin{equation*}
\begin{aligned}
&\Fc''(\e) =\y^T \Bigl(\nab^2_xB(\e\y,\cdot)\Bigr|_{v;\ v^-}^{v+v^-+q}\Bigr)\y\\
&+\y\cdot\Bigl(\frac{dq}{d\e}\nab_xB_t(\e\y,\cdot)\Bigr|_{v+v^-+q}
+\frac{dv}{d\e}\nab_xB_t(\e\y,\cdot)\Bigr|_{v}^{v+v^-+q}
\!\!+\frac{dv^-\!\!\!}{d\e}\nab_xB_t(\e\y,\cdot)\Bigr|_{v^-}^{v+v^-+q}\Bigr)\\
&+\frac{dq}{d\e}
\Bigl(\y\cdot\nab_x B_t(\e\y,\cdot)\Bigr|_{v+v^-+q}
+\frac{d(v+v^-+q)}{d\e}B_{tt}(\e\y,\cdot)\Bigr|_{v+v^-+q}\Bigr)\\
&+\frac{dv}{d\e}
\Bigl(\y\cdot\nab_x B_t(\e\y,\cdot)\Bigr|_{v}^{v+v^-+q}
\!\!+\frac{dv}{d\e}B_{tt}(\e\y,\cdot)\Bigr|_{v}^{v+v^-+q}
\!\!+\frac{d(v^-+q)}{d\e}B_{tt}(\e\y,\cdot)\Bigr|_{v+v^-+q}\Bigr)\\
&+\frac{dv^-\!\!\!}{d\e}
\Bigl(\y\cdot\nab_x B_t(\e\y,\cdot)\Bigr|_{v^-}^{v+v^-+q}
\!\!+\frac{dv^-\!\!\!}{d\e}B_{tt}(\e\y,\cdot)\Bigr|_{v^-}^{v+v^-+q}
\!\!+\frac{d(v+q)}{d\e}B_{tt}(\e\y,\cdot)\Bigr|_{v+v^-+q}\Bigr)\\
&+\frac{d^2v}{d\e^2}B_t(\e\y,\cdot)\Bigr|_{v}^{v+v^-+q}
+\frac{d^2v^-\!\!\!}{d\e^2}B_t(\e\y,\cdot)\Bigr|_{v^-}^{v+v^-+q}.
\end{aligned}
\end{equation*}


From inspection of this formula it is seen that each term in the formula
is of one of three types, which we refer to as type I, type II, or type III.
We shall invoke (\ref{dv_de_est}) to estimate them.

The only term of type I
is in the first line of this formula and is clearly $|\eta|^2O(|q|+|vv^-|)$,
by (\ref{DxB_abc}), and thus $O(1)$ by (\ref{dv_de_est}).

The terms of of type II involve, for $l=0,1$ and $k=1,2$, the quantities
\begin{equation*}
\begin{aligned}
\nab_x^l\, {\textstyle\frac{\partial^k}{\partial t^k}}\,B(\e\y,\cdot)\Bigr|_{v}^{v+v^-+q}
&=O(v^-+q)=O(e^{-c\xi^-}),\\
\nab_x^l\, {\textstyle\frac{\partial^k}{\partial t^k}}\,B(\e\y,\cdot)\Bigr|_{v^-}^{v+v^-+q}
&=O(v+q)=O(e^{-c\xi}),
\end{aligned}
\end{equation*}
which are always
multiplied by $(1+|\eta|^2)O( e^{-c\xi})$ or $(1+|\eta|^2)O( e^{-c\xi^-})$,
respectively. Thus the terms of type II are $O(1)$.

Finally the terms of of type III involve, for $l=0,1$ and $k=1,2$, the quantity
$$
\nab_x^l\, {\textstyle\frac{\partial^k}{\partial t^k}}\,B(\e\y,\cdot)\Bigr|_{v+v^-+q}
=O(1),
$$
which is always multiplied by $(1+|\eta|^2)O( e^{-c|\eta|}+e^{-c\xi}e^{-c\xi^-})$.
As above, one sees that terms of type III are bounded.
This
completes the proof of (\ref{ineqsb2}c),
and therefore (\ref{resida}).

(iii) It is sufficient to prove \eqref{residb} for $x=\xb\in\G$, as
the other case of $x\in\G^-$ is similar.
Let $\xb \in \G$ be given. Define
$s,\x^-$ and $s^-$ by the formulas
\begin{equation*}
\xb=s\ev_s=\e\x^-\ev_{r^-}+s^-\ev_{s^-}.
\end{equation*}
By (\ref{v0b}),\,(\ref{v1b}), we have $v_0(0,s)=g(\xb)-u_0(\xb)$ and $v_1(0,s)=0$;
therefore
$$
(u_0+v)\bigr|_{\xb}=u_0(\xb)+[v_0(0,s)+\e v_1(0,s)]
=g(\xb).
$$
Thus  it remains to show that $(v^-+q)\bigr|_{\xb}=O(\eps^2)$.
Indeed,  by (\ref{qt0b})\,(\ref{q1bvpb}), we have
$$
v^-+q
=[v_0^--\vo_0^-] +\e [v_1^--(\vo_1^- +\s^- \vo_{0,s^-}^-)]=O(\eps^2).
$$
In the last step here we have invoked the formulas
\begin{equation*}
\begin{aligned}
|v_0^--\vo_0^- -\e \s^- \vo_{0,s^-}^-|&=
|v_0^-(\x^-,s^-)-v_0^-(\x^-,0) -\e \s^-v_{0,s^-}^-(\x^-,0)| \\
 &= {\textstyle\half} \e^2 (\s^-)^2 |v_{0,s^-s^-}^-(\x^-,\hat s^-)|=O(\eps^2), \\
|v_1^--\vo_1^-|&=
|v_1^-(\x^-,s^-) - v_1^-(\x^-,0)|\\
   &= \e\s^-| v_{1,s^-}^-(\x^-,\hat s^-)|
   =O(\eps),
\end{aligned}
\end{equation*}
which are obtained using the exponential decay of $v_0^-$ and $v_1^-$
in $\xi^-$ and noting that
$\sigma^-=(\cot\omega)\xi^-$
on the side $\G$, where
$\omega$ is the angle at the apex.
\end{proof}

\section{Perturbed asymptotic expansion}\label{sec_beta}

The perturbed version $\b_S$
of the asymptotic expansion $\uasS$  of (\ref{uas}) is defined as follows:
\beqa
\bS(x;p)=u_0(x)+\vt(\x,s;p)+\vt^-(\x^-,s^-;p)+\qt(\y;p)+ \th p,\label{beta}
\eeqa
where a value for the positive parameter
$\theta$ and a range of values for $p$
will be chosen below.
Comparing (\ref{beta}) with (\ref{uas}), yields
$\bS(x;0)=\uasS(x) $ and,
furthermore,
an alternative equivalent representation
\begin{subequations}\label{beta_VQ}
\beq\label{beta_VQ_a}
\bS(x;p)=\uasS(x)+V(\xi,s;p)+V^-(\xi^-\!\!,s^-\!;p)+Q(\eta;p) + \th p,
\eeq
where $V=\vt-v$, $V^-=\vt^--v^-$, $Q=\qt-q$,
and therefore
\beq\label{VQ}
V=\vt_0-v_0,\qquad V^-=\vt_0^--v_0^-,\qquad Q=\qt_0-q_0.
\eeq
\end{subequations}
Note that for $V$, $V^-$ and $Q$ here, by the exponential-decay estimates for
$\frac{\pt}{\pt p}\vt_0$ and  $\frac{\pt}{\pt p}\qt_0$
from Lemma~\ref{l-v} and (\ref{q_decay}),
we have
\beq\label{VQ_decay}
\hspace{-1pt}
(1+\xi)|V|\le Cp,\quad (1+\xi^-)|V^-|\le Cp,\quad
(1+|\eta|)|Q|\le Cp e^{-c|\eta|} \le Cp.
\eeq
Furthermore, since $|\eta|\le C(\x+\x^-)$,
invoking the exponential-decay estimates for
$\frac{\pt}{\pt p}\vt_0$ and
$\frac{\pt^2}{\pt p\,\pt s}\vt_0$
from Lemma~\ref{l-v}, yields a more elaborate estimate 
\beq\label{VQ_decay_b}
(1+|\eta|e^{-c\xi^-})\,(|V|+|\textstyle\frac{\pt V}{\pt s}|)\le Cp,
\eeq
and a similar estimate involving $V^-$.

In the remainder of  this section we establish some inequalities that involve the perturbed asymptotic expansions $\b_S$.
In particular, the inequalities of Lemmas~\ref{lem_monotone} and~\ref{lem_Fbeta}
are used in \cite{KK0} to construct sub- and super-solutions
to our nonlinear boundary value problem.

\begin{lemma}\label{lem_monotone}
For the function $\bS$ of (\ref{beta}) we have
$\bS=\uasS+O(p)$. Furthermore, for some sufficiently small $\eps^*>0$,
if $p\ge0$ and $\eps\le\eps^*$,
then for all $x\in S$ we have
\beq\label{monotone}
\bS(x;-p)\le \uasS(x)-{\textstyle\frac12}\theta p,\qquad
\uasS(x)+{\textstyle\frac12}\theta p \le \bS(x;p).
\eeq
\end{lemma}

\begin{proof}
The assertion $\bS=\uasS+O(p)$ immediately follows from (\ref{beta_VQ}).
Furthermore, by (\ref{beta_VQ}), the bound for
$\bS(x;p)$ in the remaining assertion (\ref{monotone})
can be rewritten  as
\beq\label{monotone_aux}
V+V^-+Q+{\textstyle\frac12}\theta p\ge 0\qquad\mbox{for}\;\;p\ge 0.
\eeq
By (\ref{VQ_decay}), there is a sufficiently large $\bar C=\bar C(\theta)$
such that if $|\eta|\ge\bar C$, then $|Q|\le p\,C e^{-c|\eta|}\le \half\theta p$.
Combining this with $V\ge0$ and $V^-\ge 0$,
which follow from monotonicity of $\vt_0$ and $\vt_0^-$ in $p$,
established in (\ref{vt_monotone}), 
we get (\ref{monotone_aux}).
Now let $|\eta|<\bar C$. Then $|s|,\, |s^-|<\eps\bar C$.
Invoking (\ref{qdefa}), we have
$$
Q=\qt_{0}-q_0 = (\zt_{0}-z_0)-(\vto_0-\vo_0)-(\vto_0^--\vo_0^-)
=(\zt_{0}-z_0)-\mathring{V}-\mathring{V}^-,
$$
and therefore
$
V+V^-+Q=(\zt_{0}-z_0)+I+I^-,
$
where $I=V-\mathring{V}$ and $I^-=V^--\mathring{V}^-$,
and, as usual,
a small circle above a function name indicates that in the argument of the function
we have set $s=0$.
For $I$, using (\ref{VQ_decay_b}), we get
$|I|=|s\frac{\pt V}{\pt s}|\le C|s|p\le C\bar C\eps p\le \frac14 \theta p$.
Similarly, $|I^-|\le \frac14 \theta p$.
As, by Theorem~\ref{T-z0}, we also have $\zt_{0}-z_0\ge 0$,
then again $V+V^-+Q\ge -\half \theta p$.
Thus we have obtained (\ref{monotone_aux}), and therefore the bound
for $\bS(x;p)$ in (\ref{monotone}).
The bound for $\bS(x;-p)$ in (\ref{monotone}) is obtained similarly.
\end{proof}

Next, to estimate $F\b_S$, we prepare two lemmas.

\begin{lemma}
For $Q=Q(\eta;p)=\qt_0-q_0$ we have
\begin{equation}\label{t2}
\e^2\triangle Q=
B(x,\cdot)\bigr|^{\qt+\vt+\vt^-}_{\vt;\ \vt^-}
-B(x,\cdot)\bigr|^{q+v+v^-}_{v;\ v^-}
-p q_0
+O(\e^2+p^2).
\end{equation}
\end{lemma}

\begin{proof}
From (\ref{qt0a}), also using (\ref{B_tB}) and
$\qt_0=q_0+Q=q_0+O(p)$, we get
\beqa
\hspace{-0.5cm}\e^2\D Q&=&\D_\eta (\qt_0-q_0)\nonumber\\
&=&B(O,\cdot)\bigr|
^{\qt_0+\vto_0+\vto^-_0 }_{\vto_0;\ \vto_0^- }
-B(O,\cdot)\bigr|^{q_0+\vo_0+\vo^-_0}_{ \vo_0;\ \vo_0^-}
-p q_0+O(p^2).
\label{t1}
\eeqa
Recalling the definitions (\ref{VQ}), introduce the function
$$
{\cal H}(\eps,\t):=B(\e \eta,\cdot)\bigr|^{q_0+v_0+v^-_0
+\eps[q_1+v_1+v_1^-]+\t [Q+ V+V^-]}
_{ v_0+\eps v_1+\t V;\;\;  v^-_0+\eps v_1^-+\t V^-},
$$
in which we write $v_{k}=v_{k}(\xi,\e\s)$
and $v_{k}^-=v_{k}^-(\xi^-\!\!,\e\s^-)$ for $k=0,1$,
and also $V=V(\xi,\e\s)$, $V^-=V^-(\xi^-\!\!,\e\s^-)$, $Q=Q(\y)$.
The function ${\cal H}$ is defined so that, using \eqref{t1},
\begin{equation}\label{t3}
\e^2\D Q={\cal H}(0,1)-{\cal H}(0,0)-p q_0+O(p^2),
\end{equation}
and the assertion \eqref{t2} may be written as
\begin{equation}\label{t4}
\e^2\D Q={\cal H}(\eps,1)-{\cal H}(\eps,0)-p q_0+O(\eps^2+p^2).
\end{equation}
To check the formulas \eqref{t3} and \eqref{t4} we calculate
\begin{equation*}
\begin{aligned}
{\cal H}(0,0)&=B(O,\cdot)\bigr|^{q_0+\vo_0+\vo^-_0}
_{\vo_0;\;\;  \vo^-_0}, \\
{\cal H}(0,1)&=B(O,\cdot)\bigr|^{q_0+\vo_0+\vo^-_0
+ [Q+ \mathring{V}+\mathring{V}^-]}_{\vo_0+\mathring{V}\,;\;\;  \vo^-_0+\mathring{V}^-}
=B(O,\cdot)\bigr|^{\qt_0+\vto_0+\vto^-_0}_{ \vto_0;\;\;  \vto^-_0} ,\\
{\cal H}(\e,0)&=B(\e\y,\cdot)\bigr|^{q_0+v_0+v^-_0
+\e[q_1+v_1+v_1^-]}_{v_0+ \e v_1;\;\;  v^-_0+ \e v_1^-}
 =B(\e\y,\cdot)\bigr|^{q+v+v^-}_{ v;\;\; v^-}, \\
{\cal H}(\e,1)&=B(\e\y,\cdot)\bigr|^{q_0+v_0+v_0^-
+\e[q_1+v_1+v_1^-]+[Q+V+V^-]}
_{v_0+ \e v_1+V;\;\; v^-_0+ \e v_1^-+V^-}
   =B(\e\y,\cdot)\bigr|^{\qt+\vt+\vt^-}_{ \vt;\;\; \vt^-} .
\end{aligned}
\end{equation*}
Here, as usual in our notation,
a small circle above a function name indicates that in the argument of the function
we have set $s=0$; in particular,
$\mathring{V}=\vto_0-\vo_0$ and $\mathring{V}^-=\vto^-_0-\vo^-_0$.

Hence \eqref{t3} is indeed equivalent to \eqref{t1} and \eqref{t4}
is equivalent to \eqref{t2}.
To show that \eqref{t3} implies \eqref{t4} we use the mean value theorem for the
second difference and write the discrepancy between these
two formulas as
\begin{equation*}
\begin{aligned}
\Hc(\e,1)-\Hc(\e,0)-\Hc(0,1)+\Hc(0,0) = \e{\pd^2 \Hc \over \pd\e\,\pd\t}(\e^*,\t^*).
\end{aligned}
\end{equation*}
Now it suffices to show that $|{\pd^2 \Hc \over \pd\e\,\pd\t}|\le Cp$.
Then the discrepancy between the two formulas for $\e^2\triangle Q$
is bounded by $C\e|p|$, which yields
\eqref{t4}, and therefore \eqref{t2}.

To get the desired estimate for ${\pd^2 \Hc \over \pd\e\,\pd\t}$, we first evalulate
$$
{\pd\Hc \over \pd\t}(\eps,\t)=
Q\,{\cal A} +V\,{\cal B}+V^-{\cal B}^-,
$$
where
\beqann
{\cal A}&=&B_t(\e\y,\cdot)\bigr|^{q_0+v_0+v^-_0+\eps[q_1+v_1+v_1^-]+\t[Q+V+V^-]},
\\
{\cal B}&=&B_t(\e\y,\cdot)\bigr|^{q_0+v_0+v^-_0
        +\eps[q_1+v_1+v_1^-]+\t[Q+V+V^-]}_{ v_0+\eps v_1+\t V},
\\
{\cal B}^-\!\!\!\!&=&B_t(\e\y,\cdot)\bigr|^{q_0+v_0+v^-_0
        +\eps[q_1+v_1+v_1^-]+\t[Q+V+V^-]}_{v^-_0+\eps v_1^-+\t V^-}.
\eeqann
To estimate ${\pd^2\Hc \over \pd\e\pd\t}$,
we show that each of ${\pd(Q\,{\cal A}) \over \pd\e}$,
${\pd(V\,{\cal B}) \over \pd\e}$
and
${\pd(V^-{\cal B}^-) \over \pd\e}$
is $O(p)$.

For ${\pd(Q\,{\cal A}) \over \pd\e}$
we have
${\pd(Q\,{\cal A}) \over \pd\e}=Q{\pd{\cal A} \over \pd\e}$.
A calculation then shows that
$$
\begin{array}{l}
\displaystyle
{\pd{\cal A} \over \pd\e}=
\y\cdot\nab_x B_t\\
\displaystyle\quad{}+\big\{[q_1+v_1+v_1^-]
     +\s {\textstyle\frac{\pt}{\pt s}}(v_{0}+\e v_{1}+\t V)
     +\s^-\!{\textstyle\frac{\pt}{\pt s^-\!\!}}(v^-_{0}+\e v^-_{1}+\t V^-)\big\}B_{tt},
\end{array}
$$
where the terms $\nab_x B_t$ and $B_{tt}$ are computed at the point
$(\e\y,{q_0+v_0+v^-_0}+\eps[q_1+v_1+v_1^-]+\t[Q+V+V^-])$.
Recalling that $|\sigma|\le |\eta|$ and $|\s^-|\le |\eta|$, we get
$|{\pd(Q\,{\cal A}) \over \pd\e}|
\le C(1+|\eta|)|Q|\le Cp$, where we also used
(\ref{VQ_decay}).

To estimate ${\pd(V\,{\cal B}) \over \pd\e}$, another tedious calculation gives
$$
\begin{array}{l}\displaystyle
{\pd{\cal B} \over \pd\e}=
\y\cdot\nab_x B_t(\e\y,\cdot)\bigr|^{q_0+v_0+v^-_0
        +\eps[q_1+v_1+v_1^-]+\t[Q+V+V^-]}_{ v_0+\eps v_1+\t V}\\\displaystyle
{}+\big\{v_1
     +\s {\textstyle\frac{\pt}{\pt s}}(v_{0}+\e v_{1}+\t V)
     \big\}B_{tt}(\e\y,\cdot)\bigr|^{q_0+v_0+v^-_0
        +\eps[q_1+v_1+v_1^-]+\t[Q+V+V^-]}_{ v_0+\eps v_1+\t V}\\\displaystyle
{}+\big\{[q_1+v_1^-]\\\displaystyle
     \quad\;\;{}+\s^-\!{\textstyle\frac{\pt}{\pt s^-\!\!}}(v^-_{0}+\e v^-_{1}+\t V^-)\big\}
     B_{tt}
     (\e\y,\cdot)\bigr|^{q_0+v_0+v^-_0
        +\eps[q_1+v_1+v_1^-]+\t[Q+V+V^-]}.
\end{array}
$$
Now invoking Lemma~\ref{l-v} and (\ref{q_decay}),
we observe that
$|{\cal B}|\le  Ce^{-c\xi^-}$ and $|{\pd{\cal B} \over \pd\e}|
\le C(1+|\y|)e^{-c\xi^-}$.
Combining this with ${\pd V \over \pd\e}=\s {\textstyle\frac{\pt V}{\pt s}}$, where
$|\s|\le|\eta|$,
we have
$|{\pd(V\,{\cal B}) \over \pd\e}|\le C(1+|\y|)e^{-c\xi^-}(|V|+|\frac{\pt}{\pt s} V|)$
By (\ref{VQ_decay_b}), this yields the desired estimate
$|{\pd(V\,{\cal B}) \over \pd\e}|\le Cp$.
A similar argument gives $|{\pd(V^-{\cal B}^-) \over \pd\e}|\le C p$.

Thus, we have shown that each of the three components in ${\pd^2\Hc \over \pd\e\pd\t}$
is bounded by $Cp$, which completes the proof.
\end{proof}

For $F\bS$ we get the following preliminary result.

\begin{lemma}\label{L=Fbet}
For the function $\bS$ of (\ref{beta})
we have
$$
F\bS= \theta p\, b_u(x,u_0)+p \,[1+\theta\lambda(x)]\,(v_0+v_0^-+q_0)+O(\eps^2+p^2),
$$
where
$\lambda(x):=b_{uu}(x,u_0+\vartheta [v_0+v_0^-+q_0])$
with some $\vartheta=\vartheta(x)\in(0,1)$.
\end{lemma}

\begin{proof}
As from Lemma~\ref{Fuas} we have $F\uasS=O(\e^2)$,
in view of (\ref{beta_VQ}),\,(\ref{beta}) and (\ref{uas}), we get
\begin{equation*}
\begin{aligned}
F\bS &= F\bS-F\uasS+O(\e^2) \\&
 = -\e^2\D(\bS-\uasS)+b(x,\cdot)\Bigr|^{\bS}_{\uasS}+O(\e^2)\\&
 = -\e^2\D(V+V^-+Q)+B(x,\cdot)\Bigr|^{\vt+\vt^-+\qt+\theta p}_{v+v^-+q}+O(\e^2).
\end{aligned}
\end{equation*}
By (\ref{Dw0ineq}) and its analogue for $v^-$, we readily have
$$
-\e^2\D (V+V^-) = -B(x,\cdot)\Bigl|^{\vt}_{v}
\!-B(x,\cdot)\Bigl|^{\vt^-}_{v^-}
+p (v_0+v_0^-)+O(\e^2+p^2).
$$
Since (\ref{t2}) can be rewritten as
$$
-\e^2\triangle Q=
-B(x,\cdot)\bigr|^{\qt+\vt+\vt^-}_{q+v+v^-}
+B(x,\cdot)\Bigl|^{\vt}_{v}
+B(x,\cdot)\Bigl|^{\vt^-}_{v^-}
+p q_0
+O(\e^2+p^2),
$$
we now arrive at
\beqann
F\bS
 &\!\!\!\!=\!\!\!& -B(x,\cdot)\bigr|^{\qt+\vt+\vt^-}_{q+v+v^-}\!+p(v_0+v_0^-+q_0)
 +B(x,\cdot)\Bigr|^{\vt+\vt^-+\qt+\theta p}_{v+v^-+q}\!\!+O(\e^2+p^2)\\
&\!\!\!\!=\!\!\!&
B(x,\cdot)\Bigr|^{\vt+\vt^-+\qt+\theta p}_{\vt+\vt^-+\qt}+p(v_0+v_0^-+q_0)
 +O(\e^2+p^2).
 \eeqann
Note that (\ref{v}),\,(\ref{VQ_decay}) imply
that $\vt+\vt^-+\qt=v_0+v_0^-+q_0+O(\eps+p)$. Hence
\beqann
B(x,\cdot)\Bigr|^{\vt+\vt^-+\qt+\theta p}_{\vt+\vt^-+\qt}&=&
\theta p\,[B_t(x,v_0+v_0^-+q_0)+O(\eps+p)]\\
&=&\theta p\,[B_t(x,0)+\lambda(x)\,(v_0+v_0^-+q_0)+O(\eps+p)].
\eeqann
Here, by (\ref{B_tB}), one has $B_t(x,0)=b_u(x,u_0)$,
and $\lambda(x)={B_{tt}(x,\vartheta[v_0+v_0^-+q_0])}$
${{}=b_{uu}(x,u_0+\vartheta [v_0+v_0^-+q_0])}$,
as in the statement of this lemma.
Combining these formulas, we complete the proof.
\end{proof}

We are now prepared to establish our main result for $F\bS$.

\begin{lemma}\label{lem_Fbeta}
There are positive numbers $\theta$, $\eps^*$,
$p^*$ and $c_1$ such that with $\eps\le\eps^*$ and $|p| \le p^*$,
for the function $\bS$ of (\ref{beta})
one has
%
\beqann
F\b_{S} &\ge& \;\;\,\,{\textstyle\half}\theta\g^2\,p\,\,-c_1\e^2
\quad\qquad\;\;\ \!\!\rmfor p>0, 
\\
F\b_{S} &\le& -{\textstyle\half}\theta\g^2|p|+c_1\e^2
\quad\qquad\;\;\rmfor p<0. 
\eeqann
\end{lemma}

\begin{proof}
By (\ref{q_decay}), one has $|q_0| \le Ce^{-c_1|\y|}$.
Since $v_0 \ge 0$ and $v_0^- \ge 0$ it then follows that
$v_0+v_0^-+q_0 \ge -|q_0|\ge -Ce^{-c_1|\y|}$ and therefore
\begin{equation}\label{****}
v_0+v_0^-+q_0 \ge -C\i\e|\ln\e|,
\end{equation}
provided that $|\y| \ge c_1\i|\ln\e|$ and $\e^*<e\i$ so that $|\ln\e|>1$.
Furthermore, \eqref{****} also holds, with possibly a different constant $C$, when
$|\y| \le c_1\i|\ln\e|$.
Indeed, by (\ref{qdefa}),\,(\ref{T-z0-a}), we have
$v_0+v_0^-+q_0 =z_0+(v_0-\vo_0)+(v_0^--\vo_0^-) $,
where $z_0\ge0$
and, by (\ref{v}), $|v_0-\vo_0|\le C|s|$, $|v^-_0-\vo^-_0|\le C|s^-|$.
Combining these observations with
$|s|+|s^-|=\eps(|\s|+|\s^-|)\le 2\eps|\eta|$, we obtain (\ref{****})
for $|\y| \le c_1\i|\ln\e|$.
Thus we have (\ref{****}) everywhere in $S$.

Next, choose the parameter $\theta$ in the definition (\ref{beta}) of $\bS$
sufficiently small so that $0<\theta\le|\lambda(x)|^{-1}$, where
$\lambda(x)$ is from Lemma~\ref{L=Fbet}, and thus $[1+\theta\lambda]\ge0$.
Now from Lemma~\ref{L=Fbet} and \eqref{****} we obtain,
for some constants $C'$ and $C''$,
$$
F\bS\ge  \theta p \,b_u(x,u_0)-C'p\,\eps|\ln\eps|+C''(\eps^2+p^2).
$$
Since from our assumption A1 we have $b_u(x,u_0)\ge \gamma^2>0$, by choosing $\e^*$ and $p^*$
sufficiently small
we get the assertion of the lemma
 in the case $p>0$.
The case $p<0$ is similar.
\end{proof}

\section*{Conclusion}
In this note we have established four results,
Lemmas~\ref{bvps-q}, \ref{Fuas}, \ref{lem_monotone} and~\ref{lem_Fbeta},
whose proofs involve lengthy calculations.
These results
are used in \cite{KK0} to construct sub- and super-solutions
to a nonlinear boundary value problem of type (\ref{1.1}) posed in a polygonal domain.

\bibliographystyle{plain}

\end{document}